\documentclass{article}
\usepackage{amsmath}
\usepackage{amssymb}
\usepackage{amscd}
\begin{document}
\def\e#1\e{\begin{equation}#1\end{equation}}
\def\ea#1\ea{\begin{align}#1\end{align}}
\def\eq#1{{\rm(\ref{#1})}}
\newtheorem{thm}{Theorem}[section]
\newtheorem{lem}[thm]{Lemma}
\newtheorem{prop}[thm]{Proposition}
\newtheorem{conj}[thm]{Conjecture}
\newtheorem{cond}[thm]{Condition}
\newtheorem{cor}[thm]{Corollary}
\newenvironment{dfn}{\medskip\refstepcounter{thm}
\noindent{\bf Definition \thesection.\arabic{thm}\ }}{\medskip}
\newenvironment{ex}{\medskip\refstepcounter{thm}
\noindent{\bf Example \thesection.\arabic{thm}\ }}{\medskip}
\newenvironment{proof}[1][,]{\medskip\ifcat,#1
\noindent{\it Proof.\ }\else\noindent{\it Proof of #1.\ }\fi}
{\hfill$\square$\medskip}
\def\dim{\mathop{\rm dim}}
\def\Re{\mathop{\rm Re}}
\def\Im{\mathop{\rm Im}}
\def\ind{\mathop{\rm ind}}
\def\Hol{{\textstyle\mathop{\rm Hol}}}
\def\rank{\mathop{\rm rank}}
\def\diam{\mathop{\rm diam}}
\def\sign{\mathop{\rm sign}}
\def\vol{\mathop{\rm vol}}
\def\O{\mathbin{\rm O}}
\def\SO{\mathbin{\rm SO}}
\def\GL{\mathbin{\rm GL}}
\def\U{\mathbin{\rm U}}
\def\SL{\mathop{\rm SL}}
\def\SU{\mathop{\rm SU}}
\def\Sp{\mathop{\rm Sp}}
\def\Spin{\mathop{\rm Spin}}
\def\sech{{\textstyle\mathop{\rm sech}}}
\def\ge{\geqslant} 
\def\le{\leqslant} 
\def\R{\mathbin{\mathbb R}}
\def\Z{\mathbin{\mathbb Z}}
\def\C{\mathbin{\mathbb C}}
\def\CP{\mathbb{CP}}
\def\al{\alpha}
\def\be{\beta}
\def\la{\lambda}
\def\ga{\gamma}
\def\de{\delta}
\def\ep{\epsilon}
\def\ka{\kappa}
\def\th{\theta}
\def\vp{\varphi}
\def\si{\sigma}
\def\ze{\zeta}
\def\De{\Delta}
\def\La{\Lambda}
\def\Om{\Omega}
\def\Ga{\Gamma}
\def\Si{\Sigma}
\def\om{\omega}
\def\d{{\rm d}}
\def\pd{\partial}
\def\db{{\bar\partial}}
\def\ts{\textstyle}
\def\sst{\scriptscriptstyle}
\def\w{\wedge}
\def\lt{\ltimes}
\def\sm{\setminus}
\def\op{\oplus}
\def\ot{\otimes}
\def\bigot{\bigotimes}
\def\iy{\infty}
\def\ra{\rightarrow}
\def\longra{\longrightarrow}
\def\hookra{\hookrightarrow}
\def\t{\times}
\def\ha{{\textstyle\frac{1}{2}}}
\def\ti{\tilde}
\def\ovB{\,\overline{\!B}}
\def\ms#1{\vert#1\vert^2}
\def\bms#1{\bigl\vert#1\bigr\vert^2}
\def\md#1{\vert #1 \vert}
\def\bmd#1{\big\vert #1 \big\vert}
\def\nm#1{\Vert #1 \Vert}
\def\cnm#1#2{\Vert #1 \Vert_{C^{#2}}} 
\def\lnm#1#2{\Vert #1 \Vert_{L^{#2}}} 
\def\bnm#1{\bigl\Vert #1 \bigr\Vert}
\def\bcnm#1#2{\bigl\Vert #1 \bigr\Vert_{C^{#2}}} 
\def\blnm#1#2{\bigl\Vert #1 \bigr\Vert_{L^{#2}}} 
\def\an#1{\langle#1\rangle}
\def\ban#1{\bigl\langle#1\bigr\rangle}

\title{Constructing compact manifolds \\ with exceptional holonomy}
\author{Dominic Joyce \\ Lincoln College, Oxford, OX1 3DR \\
{\tt dominic.joyce@lincoln.ox.ac.uk}}
\date{March 2002}
\maketitle

\section{Introduction}
\label{s1}

In the theory of Riemannian holonomy groups, perhaps the most 
mysterious are the two exceptional cases, the holonomy group
$G_2$ in 7 dimensions and the holonomy group $\Spin(7)$ in 8 
dimensions. This is a survey paper on the construction of
examples of {\it compact} 7- and 8-manifolds with holonomy
$G_2$ and~$\Spin(7)$.

All of the material described can be found in the author's
book \cite{Joyc5}. Some, but not all, is also in the papers
\cite{Joyc1,Joyc2,Joyc3,Joyc4,Joyc6,Joyc7}. In particular,
the most complicated and powerful form of the construction
of compact manifolds with exceptional holonomy by resolving
orbifolds $T^n/\Ga$, and many of the examples, are given
only in \cite{Joyc5} and not in any published paper.

The rest of this section introduces the holonomy groups
$G_2$, $\Spin(7)$ and $\SU(m)$, and the relations between
them. Section \ref{s2} discusses constructions for compact
7-manifolds with holonomy $G_2$. Most of the section
explains how to do this by resolving the singularities
of orbifolds $T^7/\Ga$, but in \S\ref{s25} we briefly
discuss two other methods starting from Calabi--Yau 3-folds.

Section \ref{s3} explains constructions for compact 8-manifolds
with holonomy $\Spin(7)$. One way to do this is to resolve
orbifolds $T^8/\Ga$, but as this is very similar in outline
to the $G_2$ material of \S\ref{s2} we say little about it.
Instead we describe a second construction which begins with
a Calabi--Yau 4-orbifold.

\subsection{Riemannian holonomy groups}
\label{s11}

Let $M$ be a connected $n$-dimensional manifold, $g$ a Riemannian
metric on $M$, and $\nabla$ the Levi-Civita connection of $g$. Let
$x,y$ be points in $M$ joined by a smooth path $\ga$. Then {\it
parallel transport} along $\ga$ using $\nabla$ defines an isometry
between the tangent spaces $T_xM$, $T_yM$ at $x$ and~$y$.

\begin{dfn} The {\it holonomy group} $\Hol(g)$ of $g$ is the group
of isometries of $T_xM$ generated by parallel transport around
piecewise-smooth closed loops based at $x$ in $M$. We consider
$\Hol(g)$ to be a subgroup of $\O(n)$, defined up to conjugation
by elements of $\O(n)$. Then $\Hol(g)$ is independent of the base
point $x$ in~$M$.
\label{s1def1}
\end{dfn}

The classification of holonomy groups was achieved by Berger 
\cite{Berg} in~1955.

\begin{thm} Let\/ $M$ be a simply-connected, $n$-dimensional 
manifold, and\/ $g$ an irreducible, nonsymmetric Riemannian 
metric on $M$. Then either
\begin{itemize}
\setlength{\parsep}{0pt}
\setlength{\itemsep}{0pt}
\item[{\rm(i)}] $\Hol(g)=\SO(n)$,
\item[{\rm(ii)}] $n=2m$ and\/ $\Hol(g)=\SU(m)$ or\/ $\U(m)$, 
\item[{\rm(iii)}] $n=4m$ and\/ $\Hol(g)=\Sp(m)$ or\/ $\Sp(m)\Sp(1)$, 
\item[{\rm(iv)}] $n=7$ and\/ $\Hol(g)=G_2$, or
\item[{\rm(v)}] $n=8$ and\/ $\Hol(g)=\Spin(7)$.
\end{itemize}
\label{s1thm1}
\end{thm}

Now $G_2$ and $\Spin(7)$ are the exceptional cases in this 
classification, so they are called the {\it exceptional holonomy 
groups}. For some time after Berger's classification, the 
exceptional holonomy groups remained a mystery. In 1987, Bryant 
\cite{Brya} used the theory of exterior differential systems to show 
that locally there exist many metrics with these holonomy groups,
and gave some explicit, incomplete examples. Then in 1989, Bryant 
and Salamon \cite{BrSa} found explicit, {\it complete} metrics 
with holonomy $G_2$ and $\Spin(7)$ on noncompact manifolds.

In 1994-5 the author constructed the first examples of metrics
with holonomy $G_2$ and $\Spin(7)$ on {\it compact} manifolds
\cite{Joyc1,Joyc2,Joyc3}. These, and the more complicated
constructions developed later by the author \cite{Joyc4,Joyc5}
and by Kovalev \cite{Kova}, are the subject of this article.

\subsection{The holonomy group $G_2$}
\label{s12}

Let $(x_1,\dots,x_7)$ be coordinates on $\R^7$. Write
$\d{\bf x}_{ij\ldots l}$ for the exterior form $\d x_i\w\d x_j
\w\cdots\w\d x_l$ on $\R^7$. Define a metric $g_0$, a 3-form $\vp_0$
and a 4-form $*\vp_0$ on $\R^7$ by~$g_0=\d x_1^2+\cdots+\d x_7^2$,
\e
\begin{split}
\vp_0&=\d{\bf x}_{123}+\d{\bf x}_{145}+\d{\bf x}_{167}
+\d{\bf x}_{246}-\d{\bf x}_{257}-\d{\bf x}_{347}-\d{\bf x}_{356}
\;\>\text{and}\\
*\vp_0&=\d{\bf x}_{4567}+\d{\bf x}_{2367}+\d{\bf x}_{2345}
+\d{\bf x}_{1357}-\d{\bf x}_{1346}-\d{\bf x}_{1256}-\d{\bf x}_{1247}.
\end{split}
\label{s1eq1}
\e
The subgroup of $GL(7,\R)$ preserving $\vp_0$ is the {\it exceptional
Lie group} $G_2$. It also preserves $g_0,*\vp_0$ and the orientation
on $\R^7$. It is a compact, semisimple, 14-dimensional
Lie group, a subgroup of~$\SO(7)$.

A {\it $G_2$-structure} on a 7-manifold $M$ is a principal 
subbundle of the frame bundle of $M$, with structure group $G_2$. 
Each $G_2$-structure gives rise to a 3-form $\vp$ and a metric $g$ 
on $M$, such that every tangent space of $M$ admits an isomorphism 
with $\R^7$ identifying $\vp$ and $g$ with $\vp_0$ and $g_0$ 
respectively. By an abuse of notation, we will refer to $(\vp,g)$
as a $G_2$-structure.

\begin{prop} Let\/ $M$ be a $7$-manifold and\/ $(\vp,g)$ a
$G_2$-structure on $M$. Then the following are equivalent:
\begin{itemize}
\setlength{\parsep}{0pt}
\setlength{\itemsep}{0pt}
\item[{\rm(i)}] $\Hol(g)\subseteq G_2$, and\/ $\vp$ is the
induced\/ $3$-form,
\item[{\rm(ii)}] $\nabla\vp=0$ on $M$, where $\nabla$ is the
Levi-Civita connection of $g$, and
\item[{\rm(iii)}] $\d\vp=\d^*\vp=0$ on $M$.
\end{itemize}
\label{s1prop1}
\end{prop}

We call $\nabla\vp$ the {\it torsion} of the $G_2$-structure 
$(\vp,g)$, and when $\nabla\vp=0$ the $G_2$-structure is {\it 
torsion-free}. A triple $(M,\vp,g)$ is called a $G_2$-{\it manifold}
if $M$ is a 7-manifold and $(\vp,g)$ a torsion-free $G_2$-structure on
$M$. If $g$ has holonomy $\Hol(g)\subseteq G_2$, then $g$ is Ricci-flat.

\begin{thm} Let\/ $M$ be a compact\/ $7$-manifold, and suppose 
that $(\vp,g)$ is a torsion-free $G_2$-structure on $M$. 
Then $\Hol(g)=G_2$ if and only if\/ $\pi_1(M)$ is finite. In 
this case the moduli space of metrics with holonomy $G_2$ on $M$, 
up to diffeomorphisms isotopic to the identity, is a smooth 
manifold of dimension~$b^3(M)$.
\label{s1thm2}
\end{thm}

\subsection{The holonomy group $\Spin(7)$}
\label{s13}

Let $\R^8$ have coordinates $(x_1,\dots,x_8)$. Define a 4-form
$\Om_0$ on $\R^8$ by
\e
\begin{split}
\Om_0=\,&\d{\bf x}_{1234}+\d{\bf x}_{1256}+\d{\bf x}_{1278}
+\d{\bf x}_{1357}-\d{\bf x}_{1368}-\d{\bf x}_{1458}-\d{\bf x}_{1467}\\
-&\d{\bf x}_{2358}-\d{\bf x}_{2367}-\d{\bf x}_{2457}
+\d{\bf x}_{2468}+\d{\bf x}_{3456}+\d{\bf x}_{3478}+\d{\bf x}_{5678}.
\end{split}
\label{s1eq2}
\e
The subgroup of $\GL(8,\R)$ preserving $\Om_0$ is the holonomy group
$\Spin(7)$. It also preserves the orientation on $\R^8$ and the
Euclidean metric $g_0=\d x_1^2+\cdots+\d x_8^2$. It is a compact,
semisimple, 21-dimensional Lie group, a subgroup of~$\SO(8)$.

A $\Spin(7)$-structure on an 8-manifold $M$ gives rise to a 4-form 
$\Om$ and a metric $g$ on $M$, such that each tangent space of 
$M$ admits an isomorphism with $\R^8$ identifying $\Om$ and $g$ 
with $\Om_0$ and $g_0$ respectively. By an abuse of notation we 
will refer to the pair $(\Om,g)$ as a $\Spin(7)$-structure.

\begin{prop} Let\/ $M$ be an $8$-manifold and\/ $(\Om,g)$ a
$\Spin(7)$-structure on $M$. Then the following are equivalent:
\begin{itemize}
\setlength{\parsep}{0pt}
\setlength{\itemsep}{0pt}
\item[{\rm(i)}] $\Hol(g)\subseteq\Spin(7)$, and\/ $\Om$ is the 
induced\/ $4$-form,
\item[{\rm(ii)}] $\nabla\Om=0$ on $M$, where $\nabla$ is the
Levi-Civita connection of\/ $g$, and
\item[{\rm(iii)}] $\d\Om=0$ on $M$.
\end{itemize}
\label{s1prop2}
\end{prop}

We call $\nabla\Om$ the {\it torsion} of the $\Spin(7)$-structure 
$(\Om,g)$, and $(\Om,g)$ {\it torsion-free} if $\nabla\Om=0$. 
A triple $(M,\Om,g)$ is called a $\Spin(7)$-{\it manifold} if $M$
is an 8-manifold and $(\Om,g)$ a torsion-free $\Spin(7)$-structure
on $M$. If $g$ has holonomy $\Hol(g)\subseteq\Spin(7)$, then $g$
is Ricci-flat.

Here is a result on {\it compact}\/ 8-manifolds with 
holonomy~$\Spin(7)$.

\begin{thm} Let\/ $(M,\Om,g)$ be a compact\/ $\Spin(7)$-manifold.
Then $\Hol(g)=\Spin(7)$ if and only if\/ $M$ is simply-connected,
and\/ $b^3(M)+b^4_+(M)=b^2(M)+b^4_-(M)+25$. In this case the 
moduli space of metrics with holonomy $\Spin(7)$ on $M$, up to 
diffeomorphisms isotopic to the identity, is a smooth manifold of 
dimension\/~$1+b^4_-(M)$.
\label{s1thm3}
\end{thm}

\subsection{The holonomy groups $\SU(m)$}
\label{s14}

Let $\C^m\cong\R^{2m}$ have complex coordinates $(z_1,\ldots,z_m)$,
and define the metric $g_0$, K\"ahler form $\om_0$ and complex volume
form $\th_0$ on $\C^m$ by
\e
\begin{split}
g_0=\ms{\d z_1}+\cdots+\ms{\d z_m},\quad
\om_0&=\frac{i}{2}(\d z_1\w\d\bar z_1+\cdots+\d z_m\w\d\bar z_m),\\
\text{and}\quad\th_0&=\d z_1\w\cdots\w\d z_m.
\end{split}
\label{s1eq3}
\e
The subgroup of $\GL(2m,\R)$ preserving $g_0,\om_0$ and $\th_0$
is the special unitary group $\SU(m)$. Manifolds with holonomy
$\SU(m)$ are called {\it Calabi--Yau manifolds}.

Calabi--Yau manifolds are automatically Ricci-flat and K\"ahler,
with trivial canonical bundle. Conversely, any Ricci-flat K\"ahler
manifold $(M,J,g)$ with trivial canonical bundle has $\Hol(g)
\subseteq\SU(m)$. By Yau's proof of the Calabi conjecture
\cite{Yau}, we have:

\begin{thm} Let\/ $(M,J)$ be a compact complex $m$-manifold
admitting K\"ahler metrics, with trivial canonical bundle. Then
there is a unique Ricci-flat K\"ahler metric $g$ in each K\"ahler
class on $M$, and\/~$\Hol(g)\subseteq\SU(m)$.
\label{s1thm4}
\end{thm}

Using this and complex algebraic geometry one can construct many
examples of compact Calabi--Yau manifolds. The theorem also applies in
the orbifold category, yielding examples of {\it Calabi--Yau orbifolds}.

\subsection{Relations between $G_2$, $\Spin(7)$ and $\SU(m)$}
\label{s15}

Here are the inclusions between the holonomy groups $\SU(m),G_2$
and $\Spin(7)$:
\begin{equation*}
\begin{CD}
\SU(2) @>>> \SU(3) @>>> G_2 \\ @VVV @VVV @VVV \\
\SU(2)\t\SU(2) @>>> \SU(4) @>>> \Spin(7).
\end{CD}
\end{equation*}
We shall illustrate what we mean by this using the inclusion
$\SU(3)\hookra G_2$. As $\SU(3)$ acts on $\C^3$, it also acts
on $\R\op\C^3\cong\R^7$, taking the $\SU(3)$-action on $\R$ to
be trivial. Thus we embed $\SU(3)$ as a subgroup of $\GL(7,\R)$.
It turns out that $\SU(3)$ is a subgroup of the subgroup $G_2$
of $\GL(7,\R)$ defined in~\S\ref{s12}.

Here is a way to see this in terms of differential forms. Identify
$\R\op\C^3$ with $\R^7$ in the obvious way in coordinates, so that
$\bigl(x_1,(x_2+ix_3,x_4+ix_5,x_6+ix_7)\bigr)$ in $\R\op\C^3$ is
identified with $(x_1,\ldots,x_7)$ in $\R^7$. Then $\vp_0=\d x_1
\w\om_0+\Re\th_0$, where $\vp_0$ is defined in \eq{s1eq1} and
$\om_0,\th_0$ in \eq{s1eq3}. Since $\SU(3)$ preserves $\om_0$
and $\th_0$, the action of $\SU(3)$ on $\R^7$ preserves $\vp_0$,
and so~$\SU(3)\subset G_2$.

It follows that if $(M,J,h)$ is Calabi--Yau 3-fold, then $\R\t M$ and
${\cal S}^1\t M$ have torsion-free $G_2$-structures, that is, are
$G_2$-manifolds.

\begin{prop} Let\/ $(M,J,h)$ be a Calabi--Yau $3$-fold, with K\"ahler
form $\om$ and complex volume form $\th$. Let\/ $x$ be a coordinate
on $\R$ or ${\cal S}^1$. Define a metric $g=\d x^2+h$ and a $3$-form
$\vp=\d x\w\om+\Re\th$ on $\R\t M$ or ${\cal S}^1\t M$. Then
$(g,\vp)$ is a torsion-free $G_2$-structure on $\R\t M$ or
${\cal S}^1\t M$, and\/~$*\vp=\ha\om\w\om-\d x\w\Im\th$.
\label{s1prop3}
\end{prop}

Similarly, the inclusions $\SU(2)\hookra G_2$ and
$\SU(4)\hookra\Spin(7)$ give:

\begin{prop} Let\/ $(M,J,h)$ be a Calabi--Yau $2$-fold, with K\"ahler
form $\om$ and complex volume form $\th$. Let\/ $(x_1,x_2,x_3)$ be
coordinates on $\R^3$ or $T^3$. Define a metric $g=\d x_1^2+\d x_2^2+
\d x_3^2+h$ and a $3$-form
\e
\vp=\d x_1\w\d x_2\w\d x_3+\d x_1\w\om+\d x_2\w\Re\th-\d x_3\w\Im\th
\label{s1eq4}
\e
on $\R^3\t M$ or $T^3\t M$. Then $(\vp,g)$ is a torsion-free
$G_2$-structure on $\R^3\t M$ or $T^3\t M$, and
\e
*\vp=\ha\om\w\om+\d x_2\w\d x_3\w\om-\d x_1\w\d x_3\w\Re\th
-\d x_1\w\d x_2\w\Im\th.
\label{s1eq5}
\e
\label{s1prop4}
\end{prop}

\begin{prop} Let\/ $(M,J,g)$ be a Calabi--Yau $4$-fold, with K\"ahler
form $\om$ and complex volume form $\th$. Define a $4$-form $\Om$ on
$M$ by $\Om=\ha\om\w\om+\Re\th$. Then $(\Om,g)$ is a
torsion-free $\Spin(7)$-structure on~$M$.
\label{s1prop5}
\end{prop}

\section{Constructing $G_2$-manifolds from orbifolds $T^7/\Ga$\!\!}
\label{s2}

We now explain the method used in \cite{Joyc1,Joyc2} and
\cite[\S 11--\S 12]{Joyc5} to construct examples of compact
7-manifolds with holonomy $G_2$. It is based on the {\it Kummer
construction} for Calabi--Yau metrics on the $K3$ surface, and
may be divided into four steps.

\begin{list}{}{\setlength{\leftmargin}{35pt}
\setlength{\labelwidth}{35pt}}
\item[Step 1.] Let $T^7$ be the 7-torus and $(\vp_0,g_0)$ a flat
$G_2$-structure on $T^7$. Choose a finite group $\Ga$ of isometries
of $T^7$ preserving $(\vp_0,g_0)$. Then the quotient $T^7/\Ga$ is
a singular, compact 7-manifold, an {\it orbifold}.

\item[Step 2.] For certain special groups $\Ga$ there is a 
method to resolve the singularities of $T^7/\Ga$ in a natural 
way, using complex geometry. We get a nonsingular, compact 
7-manifold $M$, together with a map $\pi:M\ra T^7/\Ga$, 
the resolving map.

\item[Step 3.] On $M$, we explicitly write down a 1-parameter 
family of $G_2$-structures $(\vp_t,g_t)$ depending on $t\in(0,\ep)$.
They are not torsion-free, but have small torsion when $t$ is small.
As $t\ra 0$, the $G_2$-structure $(\vp_t,g_t)$ converges to the
singular $G_2$-structure~$\pi^*(\vp_0,g_0)$.

\item[Step 4.] We prove using analysis that for sufficiently
small $t$, the $G_2$-structure $(\vp_t,g_t)$ on $M$, with
small torsion, can be deformed to a $G_2$-structure
$(\ti\vp_t,\ti g_t)$, with zero torsion. Finally, we show
that $\ti g_t$ is a metric with holonomy $G_2$ on the compact
7-manifold~$M$.
\end{list}

We will now explain each step in greater detail.

\subsection{Step 1: Choosing an orbifold}
\label{s21}

Let $(\vp_0,g_0)$ be the Euclidean $G_2$-structure on $\R^7$
defined in \S\ref{s12}. Suppose $\La$ is a {\it lattice} in
$\R^7$, that is, a discrete additive subgroup isomorphic to
$\Z^7$. Then $\R^7/\La$ is the torus $T^7$, and $(\vp_0,g_0)$
pushes down to a torsion-free $G_2$-structure on $T^7$. We
must choose a finite group $\Ga$ acting on $T^7$ preserving
$(\vp_0,g_0)$. That is, the elements of $\Ga$ are the push-forwards
to $T^7/\La$ of affine transformations of $\R^7$ which fix
$(\vp_0,g_0)$, and take $\La$ to itself under conjugation.

Here is an example of a suitable group $\Ga$, taken from
\cite[\S 12.2]{Joyc5}.

\begin{ex} Let $(x_1,\dots,x_7)$ be coordinates
on $T^7=\R^7/\Z^7$, where $x_i\in\R/\Z$. Let $(\vp_0,g_0)$ be
the flat $G_2$-structure on $T^7$ defined by \eq{s1eq1}. Let
$\al,\be$ and $\ga$ be the involutions of $T^7$ defined by
\ea
\al:(x_1,\dots,x_7)&\mapsto(x_1,x_2,x_3,-x_4,-x_5,-x_6,-x_7),
\label{s2eq1}\\
\be:(x_1,\dots,x_7)&\mapsto(x_1,-x_2,-x_3,x_4,x_5,\ha-x_6,-x_7),
\label{s2eq2}\\
\ga:(x_1,\dots,x_7)&\mapsto\bigl(-x_1,x_2,-x_3,x_4,\ha-x_5,x_6,\ha-x_7).
\label{s2eq3}
\ea
By inspection, $\al,\be$ and $\ga$ preserve $(\vp_0,g_0)$,
because of the careful choice of exactly which signs to change. 
Also, $\al^2=\be^2=\ga^2=1$, and $\al,\be$ and  $\ga$ commute.
Thus they generate a group $\Ga=\an{\al,\be,\ga}\cong\Z_2^3$ of
isometries of $T^7$ preserving the flat $G_2$-structure~$(\vp_0,g_0)$.
\label{s2ex}
\end{ex}

Having chosen a lattice $\La$ and finite group $\Ga$, the quotient
$T^7/\Ga$ is an {\it orbifold}, a singular manifold with only quotient
singularities. The singularities of $T^7/\Ga$ come from the fixed points
of non-identity elements of $\Ga$. We now describe the singularities
in our example.

\begin{lem} In Example $\ref{s2ex}$, $\be\ga,\ga\al,\al\be$ and\/
$\al\be\ga$ have no fixed points on\/ $T^7$. The fixed points of\/
$\al,\be,\ga$ are each\/ $16$ copies of\/ $T^3$. The singular set\/
$S$ of\/ $T^7/\Ga$ is a disjoint union of\/ $12$ copies of\/ $T^3$, 
$4$ copies from each of\/ $\al,\be,\ga$. Each component of\/ $S$
is a singularity modelled on that of\/~$T^3\t\C^2/\{\pm1\}$.
\end{lem}

The most important consideration in choosing $\Ga$ is that
we should be able to resolve the singularities of $T^7/\Ga$
within holonomy $G_2$. We will explain how to do this next.

\subsection{Step 2: Resolving the singularities}
\label{s22}

Our goal is to resolve the singular set $S$ of $T^7/\Ga$ to get 
a compact 7-manifold $M$ with holonomy $G_2$. How can we do this?
In general we cannot, because we have no idea of how to resolve
general orbifold singularities with holonomy $G_2$. However,
suppose we can arrange that every connected component of $S$ is
locally isomorphic to either
\begin{itemize}
\setlength{\parsep}{0pt}
\setlength{\itemsep}{0pt}
\item[(a)] $T^3\t\C^2/G$, for $G$ a finite subgroup of $\SU(2)$, or
\item[(b)] ${\cal S}^1\t\C^3/G$, for $G$ a finite 
subgroup of $\SU(3)$ acting freely on~$\C^3\sm\{0\}$.
\end{itemize}

One can use complex algebraic geometry to find a {\it crepant
resolution} $X$ of $\C^2/G$ or $Y$ of $\C^3/G$. Then $T^3\t X$
or ${\cal S}^1\t Y$ gives a local model for how to resolve the
corresponding component of $S$ in $T^7/\Ga$. Thus we construct
a nonsingular, compact 7-manifold $M$ by using the patches
$T^3\t X$ or ${\cal S}^1\t Y$ to repair the singularities of
$T^7/\Ga$. In the case of Example \ref{s2ex}, this means gluing
12 copies of $T^3\t X$ into $T^7/\Ga$, where $X$ is the blow-up
of $\C^2/\{\pm1\}$ at its singular point.

Now the point of using crepant resolutions is this. In both case
(a) and (b), there exists a Calabi--Yau metric on $X$ or $Y$
which is asymptotic to the flat Euclidean metric on $\C^2/G$
or $\C^3/G$. Such metrics are called {\it Asymptotically Locally
Euclidean (ALE)}. In case (a), the ALE Calabi--Yau metrics were
classified by Kronheimer \cite{Kron1,Kron2}, and exist for all
finite $G\subset\SU(2)$. In case (b), crepant resolutions of
$\C^3/G$ exist for all finite $G\subset\SU(3)$ by Roan \cite{Roan},
and the author \cite{Joyc6}, \cite[\S 8]{Joyc5} proved that they
carry ALE Calabi--Yau, using a noncompact version of the Calabi
Conjecture.

By Propositions \ref{s1prop3} and \ref{s1prop4}, we can use
the Calabi--Yau metrics on $X$ or $Y$ to construct a torsion-free
$G_2$-structure on $T^3\t X$ or ${\cal S}^1\t Y$. This gives a
local model for how to resolve the singularity $T^3\t\C^2/G$ or
${\cal S}^1\t\C^3/G$ with holonomy $G_2$. So, this method gives
not only a way to smooth out the singularities of $T^7/\Ga$ as
a manifold, but also a family of torsion-free $G_2$-structures
on the resolution which show how to smooth out the singularities
of the $G_2$-structure.

The requirement above that $S$ be divided into connected components
of the form (a) and (b) is in fact unnecessarily restrictive. There
is a more complicated and powerful method, described in
\cite[\S 11--\S 12]{Joyc5}, for resolving singularities of a more
general kind. We require only that the singularities should
{\it locally} be of the form $\R^3\t\C^2/G$ or $\R\t\C^3/G$, for $G$
a finite subgroup of $\SU(2)$ or $\SU(3)$, and when $G\subset\SU(3)$
we do {\it not} require that $G$ act freely on~$\C^3\sm\{0\}$.

If $X$ is a crepant resolution of $\C^3/G$, where $G$ does not act
freely on $\C^3\sm\{0\}$, then the author shows \cite[\S 9]{Joyc5},
\cite{Joyc7} that $X$ carries a family of Calabi--Yau metrics
satisfying a complicated asymptotic condition at infinity, called
{\it Quasi-ALE} metrics. These yield the local models necessary to
resolve singularities locally of the form $\R\t\C^3/G$ with holonomy
$G_2$. Using this method we can resolve many orbifolds $T^7/\Ga$,
and prove the existence of large numbers of compact 7-manifolds
with holonomy~$G_2$.

\subsection{Step 3: Finding $G_2$-structures with small torsion}
\label{s23}

For each resolution $X$ of $\C^2/G$ in case (a), and $Y$ of 
$\C^3/G$ in case (b) above, we can find a 1-parameter family 
$\{h_t:t>0\}$ of metrics with the properties
\begin{itemize}
\setlength{\parsep}{0pt}
\setlength{\itemsep}{0pt}
\item[(a)] $h_t$ is a K\"ahler metric on $X$ with $\Hol(h_t)=\SU(2)$. 
Its injectivity radius satisfies $\de(h_t)=O(t)$, its Riemann 
curvature satisfies $\bcnm{R(h_t)}{0}=O(t^{-2})$, and 
$h_t=h+O(t^4r^{-4})$ for large $r$, where $h$ is the Euclidean 
metric on $\C^2/G$, and $r$ the distance from the origin. 
\item[(b)] $h_t$ is K\"ahler on $Y$ with $\Hol(h_t)=\SU(3)$,
where $\de(h_t)=O(t)$, $\bcnm{R(h_t)}{0}=O(t^{-2})$, and
$h_t=h+O(t^6r^{-6})$ for large~$r$.
\end{itemize}
In fact we can choose $h_t$ to be isometric to $t^2h_1$, and the
properties above are easy to prove.

Suppose one of the components of the singular set $S$ of $T^7/\Ga$
is locally modelled on $T^3\t\C^2/G$. Then $T^3$ has a natural flat
metric $h_{T^3}$. Let $X$ be the crepant resolution of $\C^2/G$ and
let $\{h_t:t>0\}$ satisfy property (a). Then Proposition \ref{s1prop4}
gives a 1-parameter family of torsion-free $G_2$-structures
$(\hat\vp_t,\hat g_t)$ on $T^3\t X$ with $\hat g_t=h_{T^3}+h_t$.
Similarly, if a component of $S$ is modelled on ${\cal S}^1\t\C^3/G$,
using Proposition \ref{s1prop3} we get a family of torsion-free
$G_2$-structures $(\hat\vp_t,\hat g_t)$ on~${\cal S}^1\t Y$.

The idea is to make a $G_2$-structure $(\vp_t,g_t)$ on $M$ by gluing 
together the torsion-free $G_2$-structures $(\hat\vp_t,\hat g_t)$ on 
the patches $T^3\t X$ and ${\cal S}^1\t Y$, and $(\vp_0,g_0)$ 
on $T^7/\Ga$. The gluing is done using a partition of unity. 
Naturally, the first derivative of the partition of unity introduces 
`errors', so that $(\vp_t,g_t)$ is not torsion-free. The size of the 
torsion $\nabla\vp_t$ depends on the difference $\hat\vp_t-\vp_0$ in 
the region where the partition of unity changes. On the patches 
$T^3\t X$, since $h_t-h=O(t^4r^{-4})$ and the partition of 
unity has nonzero derivative when $r=O(1)$, we find that 
$\nabla\vp_t=O(t^4)$. Similarly $\nabla\vp_t=O(t^6)$ on the patches
${\cal S}^1\t Y$, and so $\nabla\vp_t=O(t^4)$ on~$M$.

For small $t$, the dominant contributions to the injectivity
radius $\de(g_t)$ and Riemann curvature $R(g_t)$ are made
by those of the metrics $h_t$ on $X$ and $Y$, so we expect
$\de(g_t)=O(t)$ and $\bcnm{R(g_t)}{0}=O(t^{-2})$ by properties
(a) and (b) above. In this way we prove the following result
\cite[Th.~11.5.7]{Joyc5}, which gives the estimates on
$(\vp_t,g_t)$ that we need.

\begin{thm} On the compact\/ $7$-manifold $M$ described above,
and on many other $7$-manifolds constructed in a similar fashion,
one can write down the following data explicitly in coordinates:
\begin{itemize}
\setlength{\parsep}{0pt}
\setlength{\itemsep}{0pt}
\item Positive constants $A_1,A_2,A_3$ and\/ $\ep$,
\item A $G_2$-structure $(\vp_t,g_t)$ on $M$ with $\d\vp_t=0$
for each $t\in(0,\ep)$, and
\item A $3$-form $\psi_t$ on $M$ with\/ $\d^*\psi_t=\d^*\vp_t$
for each\/~$t\in(0,\ep)$.
\end{itemize}

These satisfy three conditions:
\begin{itemize}
\setlength{\parsep}{0pt}
\setlength{\itemsep}{0pt}
\item[{\rm(i)}] $\lnm{\psi_t}{2}\le A_1t^4$, $\cnm{\psi_t}{0}
\le A_1t^3$ and\/~$\lnm{\d^*\psi_t}{14}\le A_1t^{16/7}$,
\item[{\rm(ii)}] the injectivity radius $\de(g_t)$ satisfies
$\de(g_t)\ge A_2t$,
\item[{\rm(iii)}] the Riemann curvature $R(g_t)$ of\/ $g_t$
satisfies~$\bcnm{R(g_t)}{0}\le A_3t^{-2}$.
\end{itemize}
Here the operator $\d^*$ and the norms $\lnm{\,.\,}{2}$,
$\lnm{\,.\,}{14}$ and\/ $\cnm{\,.\,}{0}$ depend on~$g_t$.
\label{s2thm1}
\end{thm}

Here one should regard $\psi_t$ as a {\it first integral} of 
the torsion $\nabla\vp_t$ of $(\vp_t,g_t)$. Thus the norms 
$\lnm{\psi_t}{2}$, $\cnm{\psi_t}{0}$  and $\lnm{\d^*\psi_t}{14}$
are measures of $\nabla\vp_t$. So parts (i)--(iii) say that the 
torsion $\nabla\vp_t$ must be small compared to the injectivity 
radius and Riemann curvature of~$(M,g_t)$.

\subsection{Step 4: Deforming to a torsion-free $G_2$-structure}
\label{s24}

We prove the following analysis result.

\begin{thm} Let\/ $A_1,A_2,A_3$ be positive constants. Then 
there exist positive constants $\ka,K$ such that whenever
$0<t\le\ka$, the following is true.

Let\/ $M$ be a compact\/ $7$-manifold, and\/ $(\vp,g)$ a 
$G_2$-structure on $M$ with\/ $\d\vp\!=\!0$. Suppose $\psi$ is 
a smooth\/ $3$-form on\/ $M$ with\/ $\d^*\psi=\d^*\vp$, and
\begin{itemize}
\setlength{\parsep}{0pt}
\setlength{\itemsep}{0pt}
\item[{\rm(i)}] $\lnm{\psi}{2}\le A_1 t^4$, 
$\cnm{\psi}{0}\le A_1 t^{1/2}$
and\/~$\lnm{\d^*\psi}{14}\le A_1$,
\item[{\rm(ii)}] the injectivity radius $\de(g)$ satisfies
$\de(g)\ge A_2 t$, and
\item[{\rm(iii)}] the Riemann curvature $R(g)$ 
satisfies~$\bcnm{R(g)}{0}\le A_3 t^{-2}$.
\end{itemize}
Then there exists a smooth, torsion-free $G_2$-structure 
$(\ti\vp,\ti g)$ on $M$ with\/ $\cnm{\ti\vp-\vp}{0}\le Kt^{1/2}$.
\label{s2thm2}
\end{thm}

Basically, this result says that if $(\vp,g)$ is a $G_2$-structure 
on $M$, and the torsion $\nabla\vp$ is sufficiently small, then we 
can deform to a nearby $G_2$-structure $(\ti\vp,\ti g)$ that 
is torsion-free. Here is a sketch of the proof of Theorem \ref{s2thm2}, 
ignoring several technical points. The proof is that given in 
\cite[\S 11.6--\S 11.8]{Joyc5}, which is an improved version of
the proof in~\cite{Joyc1}.

We have a 3-form $\vp$ with $\d\vp=0$ and $\d^*\vp=\d^*\psi$ for small 
$\psi$, and we wish to construct a nearby 3-form $\ti\vp$ with 
$\d\ti\vp=0$ and $\ti \d^*\ti\vp=0$. Set $\ti\vp=\vp+\d\eta$, 
where $\eta$ is a small 2-form. Then $\eta$ must satisfy a nonlinear 
p.d.e., which we write as
\e
\d^*\d\eta=-\d^*\psi+\d^*F(\d\eta),
\label{s2eq4}
\e
where $F$ is nonlinear, satisfying~$F(\d\eta)=O\bigl(\ms{\d\eta}\bigr)$.

We solve \eq{s2eq4} by iteration, introducing a sequence
$\{\eta_j\}_{j=0}^\infty$ with $\eta_0=0$, satisfying the
inductive equations
\e
\d^*\d\eta_{j+1}=-\d^*\psi+\d^*F(\d\eta_j),\qquad\qquad
\d^*\eta_{j+1}=0.
\label{s2eq5}
\e
If such a sequence exists and converges to $\eta$, then taking
the limit in \eq{s2eq5} shows that $\eta$ satisfies 
\eq{s2eq4}, giving us the solution we want. 

The key to proving this is an {\it inductive estimate} on the
sequence $\{\eta_j\}_{j=0}^\infty$. The inductive estimate we use 
has three ingredients, the equations
\ea
\lnm{\d\eta_{j+1}}{2}&\le\lnm{\psi}{2}
+C_1\lnm{\d\eta_j}{2}\cnm{\d\eta_j}{0},
\label{s2eq6}\\
\lnm{\nabla\d\eta_{j+1}}{14}&\le
C_2\bigl(\lnm{\d^*\psi}{14}+
\lnm{\nabla\d\eta_j}{14}\cnm{\d\eta_j}{0}
+t^{-4}\lnm{\d\eta_{j+1}}{2}\bigr),
\label{s2eq7}\\
\cnm{\d\eta_j}{0}&\le C_3\bigl(t^{1/2}\lnm{\nabla\d\eta_j}{14}
+t^{-7/2}\lnm{\d\eta_j}{2}\bigr).
\label{s2eq8}
\ea
Here $C_1,C_2,C_3$ are positive constants independent of $t$. Equation
\eq{s2eq6} is obtained from \eq{s2eq5} by taking the $L^2$-inner 
product with $\eta_{j+1}$ and integrating by parts. Using the fact that
$\d^*\vp=\d^*\psi$ and $\lnm{\psi}{2}=O(t^4)$, $\md{\psi}=O(t^{1/2})$
we get a powerful estimate of the $L^2$-norm of~$\d\eta_{j+1}$.

Equation \eq{s2eq7} is derived from an {\it elliptic regularity 
estimate} for the operator $\d+\d^*$ acting on 3-forms on $M$. Equation
\eq{s2eq8} follows from the {\it Sobolev embedding theorem}, 
since $L^{14}_1(M)$ embeds in $C^0(M)$. Both \eq{s2eq7} and
\eq{s2eq8} are proved on small balls of radius $O(t)$ in $M$,
using parts (ii) and (iii) of Theorem \ref{s2thm1}, and this
is where the powers of $t$ come from.

Using \eq{s2eq6}-\eq{s2eq8} and part (i) of Theorem \ref{s2thm1}
we show that if 
\e
\lnm{\d\eta_j}{2}\le C_4t^4,\;\>
\lnm{\nabla\d\eta_j}{14}\le C_5,\;\>\text{and}\;\>
\cnm{\d\eta_j}{0}\le Kt^{1/2},
\label{s2eq9}
\e
where $C_4,C_5$ and $K$ are positive constants depending on
$C_1,C_2,C_3$ and $A_1$, and if $t$ is sufficiently small, then 
the same inequalities \eq{s2eq9} apply to $\d\eta_{j+1}$. Since 
$\eta_0=0$, by induction \eq{s2eq9} applies for all $j$ and the 
sequence $\{\d\eta_j\}_{j=0}^\infty$ is bounded in the Banach space 
$L^{14}_1(\Lambda^3T^*M)$. One can then use standard techniques in 
analysis to prove that this sequence converges to a smooth 
limit $\d\eta$. This concludes the proof of Theorem~\ref{s2thm2}.

\begin{figure}[htb]
{\caption{Betti numbers $(b^2,b^3)$ of compact
$G_2$-manifolds}\label{g2betti}}
{\setlength{\unitlength}{1.5pt}
\begin{picture}(215,110)(-7,-10)
\linethickness{0.25mm}
\put(0,0){\line(1,0){220}}
\put(0,0){\line(0,1){90}}
\multiput(0,15)(0,15){5}{\multiput(1,0)(1,0){220}{\vrule width .25pt
height.125pt depth.125pt}}
\multiput(10,0)(10,0){21}{\multiput(0,3)(0,3){30}{\vrule width .25pt
height.25pt depth0pt}}
\def\rla#1{\hbox to 0pt{\footnotesize\hss #1}}
\def\cla#1{\hbox to 0pt{\footnotesize\hss #1\hss}}
\put(-2,-2){\rla{0}}\put(-2,13){\rla{5}}\put(-2,28){\rla{10}}
\put(-2,43){\rla{15}}\put(-2,58){\rla{20}}\put(-2,73){\rla{25}}
\put(0,-6){\cla{0}}\put(10,-6){\cla{10}}\put(20,-6){\cla{20}}
\put(30,-6){\cla{30}}\put(40,-6){\cla{40}}\put(50,-6){\cla{50}}
\put(60,-6){\cla{60}}\put(70,-6){\cla{70}}\put(80,-6){\cla{80}}
\put(90,-6){\cla{90}}\put(100,-6){\cla{100}}\put(110,-6){\cla{110}}
\put(120,-6){\cla{120}}\put(130,-6){\cla{130}}\put(140,-6){\cla{140}}
\put(150,-6){\cla{150}}\put(160,-6){\cla{160}}\put(170,-6){\cla{170}}
\put(180,-6){\cla{180}}\put(190,-6){\cla{190}}\put(200,-6){\cla{200}}
\put(210,-6){\cla{210}}
\put(-4,95){$b^2(M)$}\put(200,6){$b^3(M)$}
\def\q#1#2{\count1=#1 \multiply\count1 by 3 \put(#2,\count1){\circle*{1.6}}}
\q{8}{47}\q{9}{46}\q{10}{45}
\q{11}{44}\q{12}{43}\q{13}{42}
\q{14}{41}\q{15}{40}\q{16}{39}
\q{4}{35}\q{5}{34}
\q{6}{33}\q{7}{32}\q{8}{31}
\q{9}{30}\q{10}{29}\q{11}{28}
\q{12}{27}\q{11}{36}
\q{5}{13}\q{3}{11}
\q{4}{17}\q{2}{10}
\q{5}{18}\q{6}{17}\q{3}{6}
\q{2}{7}\q{8}{7}\q{3}{4} 
\q{4}{99}\q{4}{95}\q{4}{91}\q{4}{87}\q{4}{83}\q{6}{89} 
\q{4}{79}\q{6}{85}\q{4}{75}\q{6}{81}\q{4}{71}\q{6}{77} 
\q{8}{83}\q{4}{67}\q{6}{73}\q{8}{79}\q{4}{63}\q{6}{69} 
\q{8}{75}\q{4}{59}\q{6}{65}\q{8}{71}\q{10}{77}\q{4}{55} 
\q{6}{61}\q{8}{67}\q{10}{73}\q{4}{51}\q{6}{57}\q{8}{63} 
\q{10}{69}\q{12}{75}\q{4}{47}\q{6}{53}\q{8}{59}\q{10}{65} 
\q{12}{71}\q{4}{43}\q{6}{49}\q{8}{55}\q{10}{61}\q{12}{67} 
\q{14}{73}\q{4}{39}\q{6}{45}\q{8}{51}\q{10}{57}\q{12}{63} 
\q{14}{69}\q{16}{75}\q{4}{35}\q{6}{41}\q{8}{47}\q{10}{53} 
\q{12}{59}\q{14}{65}\q{16}{71}\q{18}{77}\q{2}{157}\q{3}{128} 
\q{4}{99}\q{2}{153}\q{3}{124}\q{4}{95}\q{2}{149}\q{3}{120} 
\q{4}{91}\q{2}{145}\q{3}{116}\q{4}{87}\q{2}{141}\q{4}{147} 
\q{3}{112}\q{5}{118}\q{4}{83}\q{6}{89}\q{2}{137}\q{4}{143} 
\q{3}{108}\q{5}{114}\q{4}{79}\q{6}{85}\q{2}{133}\q{4}{139} 
\q{3}{104}\q{5}{110}\q{4}{75}\q{6}{81}\q{2}{129}\q{4}{135} 
\q{6}{141}\q{3}{100}\q{5}{106}\q{7}{112}\q{4}{71}\q{6}{77} 
\q{8}{83}\q{2}{125}\q{4}{131}\q{6}{137}\q{3}{96}\q{5}{102} 
\q{7}{108}\q{4}{67}\q{6}{73}\q{8}{79}\q{2}{121}\q{4}{127} 
\q{6}{133}\q{3}{92}\q{5}{98}\q{7}{104}\q{4}{63}\q{6}{69} 
\q{8}{75}\q{4}{123}\q{6}{129}\q{8}{135}\q{5}{94}\q{7}{100} 
\q{9}{106}\q{6}{65}\q{8}{71}\q{10}{77}\q{6}{125}\q{8}{131} 
\q{7}{96}\q{9}{102}\q{8}{67}\q{10}{73}\q{8}{127}\q{10}{133} 
\q{9}{98}\q{11}{104}\q{10}{69}\q{12}{75}\q{10}{129}\q{11}{100} 
\q{12}{71}\q{12}{131}\q{13}{102}\q{14}{73}\q{14}{133}\q{15}{104} 
\q{16}{75}\q{16}{135}\q{17}{106}\q{18}{77}\q{0}{215}\q{1}{186} 
\q{2}{157}\q{3}{128}\q{4}{99}
\q{5}{74}\q{5}{70}\q{5}{66}\q{5}{62}\q{5}{58}\q{7}{64}
\q{6}{65}\q{6}{61}\q{6}{57}\q{6}{53}\q{6}{49}\q{8}{55}
\q{7}{56}\q{7}{52}\q{7}{48}\q{7}{44}\q{7}{40}\q{9}{46}
\q{9}{46}\q{8}{47}\q{8}{43}\q{8}{39}\q{8}{35}\q{8}{31}
\q{10}{37}\q{10}{53}\q{10}{49}\q{10}{45}\q{10}{41}
\q{10}{37}\q{12}{43}\q{12}{55}\q{12}{51}\q{12}{47}
\q{12}{43}\q{11}{44}\q{11}{40}\q{11}{36}\q{11}{32}
\q{11}{28}\q{13}{34}\q{13}{46}\q{13}{42}\q{13}{38}
\q{13}{34}\q{15}{40}\q{14}{41}\q{14}{37}\q{14}{33}
\q{14}{29}\q{14}{25}\q{16}{31}\q{16}{43}\q{16}{39}
\q{16}{35}\q{16}{31}\q{18}{45}\q{18}{41}\q{18}{37}
\q{20}{43}\q{17}{38}\q{17}{34}\q{17}{30}\q{17}{26}
\q{19}{28}\q{19}{40}\q{19}{36}\q{19}{32}\q{19}{28}
\q{21}{42}\q{21}{38}\q{21}{34}\q{23}{44}\q{23}{40}
\q{25}{46}\q{3}{84}\q{3}{80}\q{3}{76}\q{3}{72}\q{3}{68}
\q{4}{79}\q{4}{75}\q{4}{71}\q{4}{67}\q{4}{63}\q{4}{75}
\q{4}{71}\q{4}{67}\q{4}{63}\q{4}{59}\q{5}{70}\q{5}{66}
\q{5}{62}\q{5}{58}\q{5}{54}\q{5}{66}\q{5}{62}\q{5}{58}
\q{5}{54}\q{5}{50}\q{6}{61}\q{6}{57}\q{6}{53}\q{6}{49}
\q{6}{45}\q{7}{56}\q{8}{51}\q{6}{57}\q{6}{53}\q{6}{49}
\q{6}{45}\q{7}{52}\q{7}{48}\q{7}{44}\q{7}{40}\q{8}{63}
\q{8}{59}\q{8}{55}\q{8}{51}\q{8}{47}\q{9}{58}\q{9}{54}
\q{9}{50}\q{9}{46}\q{9}{42}\q{10}{65}\q{10}{61}\q{10}{57}
\q{10}{53}\q{11}{60}\q{11}{56}\q{11}{52}\q{11}{48}\q{9}{54}
\q{9}{50}\q{9}{46}\q{10}{49}\q{10}{45}\q{10}{41}\q{11}{56}
\q{11}{52}\q{11}{48}\q{12}{51}\q{12}{47}\q{12}{43}\q{13}{50}
\q{14}{45}\q{12}{51}\q{12}{47}\q{13}{46}\q{13}{42}\q{14}{53}
\q{14}{49}\q{15}{48}\q{15}{44}\q{16}{55}\q{16}{51}\q{17}{50}
\q{17}{46}\q{18}{53}\q{19}{48}\q{15}{48}\q{16}{43}\q{17}{50}
\q{18}{45}\q{19}{52}\q{20}{47}\q{21}{54}\q{22}{49}\q{23}{56}
\q{24}{51}\q{1}{142}\q{2}{137}\q{3}{132}\q{3}{128}\q{3}{124}
\q{3}{120}\q{5}{122}\q{2}{113}\q{3}{108}\q{4}{103}\q{4}{99}
\q{4}{95}\q{4}{91}\q{6}{93}\q{3}{84}\q{4}{79}\q{5}{74}
\q{5}{70}\q{5}{66}\q{5}{62}\q{7}{64}
\q{28}{43}
\q{5}{26}\q{10}{29}\q{15}{32}\q{20}{35}
\q{10}{13}\q{5}{10}
\q{2}{19}\q{3}{19}\q{4}{19}\q{5}{19}\q{6}{19}\q{7}{19}\q{8}{19}
\q{9}{19}\q{10}{19}\q{11}{19}
\end{picture}}
\end{figure}

From Theorems \ref{s2thm1} and \ref{s2thm2} we see that the compact
7-manifold $M$ constructed in Step 2 admits torsion-free
$G_2$-structures $(\ti\vp,\ti g)$. Theorem \ref{s1thm2} then
shows that $\Hol(\ti g)=G_2$ if and only if $\pi_1(M)$ is finite.
In the example above $M$ is simply-connected, and so $\pi_1(M)=\{1\}$
and $M$ has metrics with holonomy $G_2$, as we want.

By considering different groups $\Ga$ acting on $T^7$, and also 
by finding topologically distinct resolutions $M_1,\dots,M_k$ of 
the same orbifold $T^7/\Ga$, we can construct many compact 
Riemannian 7-manifolds with holonomy $G_2$. A good number of
examples are given in \cite[\S 12]{Joyc5}. Figure \ref{g2betti}
displays the Betti numbers of compact, simply-connected 7-manifolds
with holonomy $G_2$ constructed there. There are 252 different sets
of Betti numbers.

Examples are also known \cite[\S 12.4]{Joyc5} of compact 7-manifolds
with holonomy $G_2$ with finite, nontrivial fundamental group. It
seems likely to the author that the Betti numbers given in Figure
\ref{g2betti} are only a small proportion of the Betti numbers of
all compact, simply-connected 7-manifolds with holonomy~$G_2$.

\subsection{Other constructions of compact $G_2$-manifolds}
\label{s25}

Here are two other methods, taken from \cite[\S 11.9]{Joyc5},
that may be used to construct compact 7-manifolds with
holonomy $G_2$. The first method was outlined by the author
in \cite[\S 4.3]{Joyc2}, and is being studied by the author's
student Ben Stephens.
\medskip

\noindent{\bf Method 1.} Let $(Y,J,h)$ be a Calabi--Yau 3-fold, 
with K\"ahler form $\om$ and holomorphic volume form $\th$. Suppose
$\si:Y\ra Y$ is an involution, satisfying $\si^*(h)=h$,
$\si^*(J)=-J$ and $\si^*(\th)=\bar\th$. We call
$\si$ a {\it real structure} on $Y$. Let $N$ be the fixed point 
set of $\si$ in $Y$. Then $N$ is a real 3-dimensional submanifold 
of $Y$, and is in fact a special Lagrangian 3-fold.

Let ${\cal S}^1=\R/\Z$, and define a torsion-free $G_2$-structure 
$(\vp,g)$ on ${\cal S}^1\t Y$ as in Proposition \ref{s1prop3}. 
Then $\vp=\d x\w\om+\Re\th$, where $x\in\R/\Z$ is the coordinate 
on ${\cal S}^1$. Define $\hat\si:{\cal S}^1\t Y\ra{\cal S}^1\t Y$ 
by $\hat\si\bigl((x,y)\bigr)=\bigl(-x,\si(y)\bigr)$. Then $\hat\si$
preserves $(\vp,g)$ and $\hat\si^2=1$. The fixed points of 
$\hat\si$ in ${\cal S}^1\t Y$ are $\{\Z,\ha+\Z\}\t N$. Thus 
$({\cal S}^1\t Y)/\an{\hat\si}$ is an orbifold. Its singular
set is 2 copies of $N$, and each singular point is modelled 
on~$\R^3\t\R^4/\{\pm 1\}$.

We aim to resolve $({\cal S}^1\t Y)/\an{\hat\si}$ to get a
compact 7-manifold $M$ with holonomy $G_2$. Locally, each
singular point should be resolved like $\R^3\t X$, where
$X$ is an ALE Calabi--Yau 2-fold asymptotic to $\C^2/\{\pm 1\}$.
There is a 3-dimensional family of such $X$, and we need to
choose one member of this family for each singular point in
the singular set.

Calculations by the author indicate that the data needed to do 
this is a closed, coclosed 1-form $\al$ on $N$ that is nonzero
at every point of $N$. The existence of a suitable 1-form $\al$
depends on the metric on $N$, which is the restriction of the
metric $g$ on $Y$. But $g$ comes from the solution of the Calabi
Conjecture, so we know little about it. This may make the method
difficult to apply in practice.
\medskip

The second method is studied by Alexei Kovalev \cite{Kova}, and
is based on an idea due to Simon Donaldson.
\medskip

\noindent{\bf Method 2.} Let $X$ be a projective complex 3-fold
with canonical bundle $K_X$, and $s$ a holomorphic section of 
$K_X^{-1}$ which vanishes to order 1 on a smooth divisor $D$ 
in $X$. Then $D$ has trivial canonical bundle, so $D$ is $T^4$
or $K3$. Suppose $D$ is a $K3$ surface. Define $Y=X\sm D$,
and suppose $Y$ is simply-connected.

Then $Y$ is a noncompact complex 3-fold with $K_Y$ trivial, 
and one infinite end modelled on $D\t{\cal S}^1\t[0,\iy)$.
Using a version of the proof of the Calabi Conjecture for
noncompact manifolds one constructs a complete Calabi--Yau
metric $h$ on $Y$, which is asymptotic to the product on
$D\t{\cal S}^1\t[0,\iy)$ of a Calabi--Yau metric on $D$,
and Euclidean metrics on ${\cal S}^1$ and $[0,\iy)$. We
call such metrics {\it Asymptotically Cylindrical}.

Suppose we have such a metric on $Y$. Define a torsion-free 
$G_2$-structure $(\vp,g)$ on ${\cal S}^1\t Y$ as in Proposition
\ref{s1prop3}. Then ${\cal S}^1\t Y$ is a noncompact $G_2$-manifold
with one end modelled on $D\t T^2\t[0,\iy)$, whose metric is
asymptotic to the product on $D\t T^2\t[0,\iy)$ of a Calabi--Yau
metric on $D$, and Euclidean metrics on $T^2$ and~$[0,\iy)$. 

Donaldson and Kovalev's idea is to take two such products
${\cal S}^1\t Y_1$ and ${\cal S}^1\t Y_2$ whose infinite ends 
are isomorphic in a suitable way, and glue them together to get 
a compact 7-manifold $M$ with holonomy $G_2$. The gluing process
swaps round the ${\cal S}^1$ factors. That is, the ${\cal S}^1$
factor in ${\cal S}^1\t Y_1$ is identified with the asymptotic
${\cal S}^1$ factor in $Y_2\sim D_2\t{\cal S}^1\t[0,\iy)$, and
vice versa.

\section{Compact $\Spin(7)$-manifolds from Calabi--Yau 4-orbifolds}
\label{s3}

In a very similar way to the $G_2$ case, one can construct examples
of compact 8-manifolds with holonomy $\Spin(7)$ by resolving the
singularities of torus orbifolds $T^8/\Ga$. This is done in
\cite{Joyc3} and \cite[\S 13--\S 14]{Joyc5}. In \cite[\S 14]{Joyc5},
examples are constructed which realize 181 different sets of Betti
numbers. Two compact 8-manifolds with holonomy $\Spin(7)$ and the
same Betti numbers may be distinguished by the cup products on their
cohomologies (examples of this are given in \cite[\S 3.4]{Joyc3}),
so they probably represent rather more than 181 topologically
distinct 8-manifolds.

The main differences with the $G_2$ case are, firstly, that the
technical details of the analysis are different and harder, and
secondly, that the singularities that arise are typically more
complicated and more tricky to resolve. One reason for this is that
in the $G_2$ case the singular set is made up of 1 and 3-dimensional
pieces in a 7-dimensional space, so one can often arrange for the
pieces to avoid each other, and resolve them independently.

But in the $\Spin(7)$ case the singular set is typically made up
of 4-dimensional pieces in an 8-dimensional space, so they nearly
always intersect. There are also topological constraints arising
from the $\hat A$-genus, which do not apply in the $G_2$ case.
The moral appears to be that when you increase the dimension,
things become more difficult.

Anyway, we will not discuss this further, as the principles are
very similar to the $G_2$ case above. Instead, we will discuss
an entirely different construction of compact 8-manifolds with
holonomy $\Spin(7)$ developed by the author in \cite{Joyc4} and
\cite[\S 15]{Joyc5}, a little like Method 1 of \S\ref{s25}. In
this we start from a {\it Calabi--Yau $4$-orbifold} rather than
from $T^8$. The construction can be divided into five steps.
\begin{list}{}{\setlength{\leftmargin}{40pt}
\setlength{\labelwidth}{40pt}}
\item[Step 1.] Find a compact, complex 4-orbifold $(Y,J)$
satisfying the conditions:
\begin{itemize}
\setlength{\parsep}{0pt}
\setlength{\itemsep}{0pt}
\item[(a)] $Y$ has only finitely many singular points
$p_1,\ldots,p_k$, for~$k\ge 1$.
\item[(b)] $Y$ is modelled on $\C^4/\an{i}$ near each $p_j$,
where $i$ acts on $\C^4$ by complex multiplication.
\item[(c)] There exists an antiholomorphic involution
$\si:Y\ra Y$ whose fixed point set is~$\{p_1,\ldots,p_k\}$.
\item[(d)] $Y\sm\{p_1,\ldots,p_k\}$ is simply-connected,
and~$h^{2,0}(Y)=0$.
\end{itemize}
\item[Step 2.] Choose a $\si$-invariant K\"ahler class on $Y$.
Then by Theorem \ref{s1thm4} there exists a unique $\si$-invariant
Ricci-flat K\"ahler metric $g$ in this K\"ahler class. Let $\om$
be the K\"ahler form of $g$. Let $\th$ be a holomorphic volume
form for $(Y,J,g)$. By multiplying $\th$ by ${\rm e}^{i\phi}$ if
necessary, we can arrange that~$\si^*(\th)=\bar\th$. 

Define $\Om=\ha\om\w\om+\Re\th$. Then $(\Om,g)$ is a torsion-free
$\Spin(7)$-structure on $Y$, by Proposition \ref{s1prop5}. Also,
$(\Om,g)$ is $\si$-invariant, as $\si^*(\om)=-\om$ and $\si^*(\th)
=\bar\th$. Define $Z=Y/\an{\si}$. Then $Z$ is a compact real
8-orbifold with isolated singular points $p_1,\ldots,p_k$, and
$(\Om,g)$ pushes down to a torsion-free $\Spin(7)$-structure
$(\Om,g)$ on $Z$.
\item[Step 3.] $Z$ is modelled on $\R^8/G$ near each $p_j$, where
$G$ is a certain finite subgroup of $\Spin(7)$ with $\md{G}=8$.
We can write down two explicit, topologically distinct ALE
$\Spin(7)$-manifolds $X_1,X_2$ asymptotic to $\R^8/G$. Each
carries a 1-parameter family of homothetic ALE metrics $h_t$
for $t>0$ with $\Hol(h_t)=\Z_2\lt\SU(4)\subset\Spin(7)$.

For $j=1,\ldots,k$ we choose $i_j=1$ or 2, and resolve the
singularities of $Z$ by gluing in $X_{i_j}$ at the singular
point $p_j$ for $j=1,\ldots,k$, to get a compact, nonsingular
8-manifold $M$, with projection~$\pi:M\ra Z$.

\item[Step 4.] On $M$, we explicitly write down a 1-parameter 
family of $\Spin(7)$-structures $(\Om_t,g_t)$ depending on $t\in
(0,\ep)$. They are not torsion-free, but have small torsion when
$t$ is small. As $t\ra 0$, the $\Spin(7)$-structure $(\Om_t,g_t)$
converges to the singular $\Spin(7)$-structure~$\pi^*(\Om_0,g_0)$.
\item[Step 5.] We prove using analysis that for sufficiently
small $t$, the $\Spin(7)$-structure $(\Om_t,g_t)$ on $M$, with
small torsion, can be deformed to a $\Spin(7)$-structure
$(\ti\Om_t,\ti g_t)$, with zero torsion.

It turns out that if $i_j=1$ for $j=1,\ldots,k$ we have
$\pi_1(M)\cong\Z_2$ and $\Hol(\ti g_t)=\Z_2\lt\SU(4)$, and for
the other $2^k-1$ choices of $i_1,\ldots,i_k$ we have $\pi_1(M)=\{1\}$
and $\Hol(\ti g_t)=\Spin(7)$. So $\ti g_t$ is a metric with holonomy
$\Spin(7)$ on the compact 8-manifold $M$ for~$(i_1,\ldots,i_k)\ne
(1,\ldots,1)$.
\end{list}

Once we have completed Step 1, Step 2 is immediate. Steps 4 and 5
are analogous to Steps 3 and 4 of \S\ref{s2}, and can be done using
the techniques and analytic results developed by the author for
the first $T^8/\Ga$ construction of compact $\Spin(7)$-manifolds,
\cite{Joyc3}, \cite[\S 13]{Joyc5}. So the really new material is
in Steps 1 and 3, and we will discuss only these.

\subsection{Step 1: An example}
\label{s31}

We do Step 1 using complex algebraic geometry. The problem is that
conditions (a)--(d) above are very restrictive, so it is not that
easy to find {\it any} $Y$ satisfying all four conditions. All the
examples $Y$ the author has found are constructed using {\it weighted
projective spaces}, an important class of complex orbifolds.

\begin{dfn} Let $m\ge 1$ be an integer, and $a_0,a_1,\ldots,a_m$
positive integers with highest common factor 1. Let $\C^{m+1}$ have 
complex coordinates on $(z_0,\ldots,z_m)$, and define an action of 
the complex Lie group $\C^*$ on $\C^{m+1}$ by
\begin{equation*}
(z_0,\ldots,z_m)\,{\buildrel u\over\longmapsto}
(u^{a_0}z_0,\ldots,u^{a_m}z_m),\qquad\text{for $u\in\C^*$.}
\end{equation*}
The {\it weighted projective space} $\CP^m_{a_0,\ldots,a_m}$
is $\bigl(\C^{m+1}\sm\{0\}\bigr)/\C^*$. The $\C^*$-orbit of
$(z_0,\ldots,z_m)$ is written~$[z_0,\ldots,z_m]$.
\label{s3def}
\end{dfn}

Here is the simplest example the author knows.

\begin{ex} Let $Y$ be the hypersurface of degree 12 in 
$\CP^5_{1,1,1,1,4,4}$ given by
\begin{equation*}
Y=\bigl\{[z_0,\ldots,z_5]\in\CP^5_{1,1,1,1,4,4}:z_0^{12}+
z_1^{12}+z_2^{12}+z_3^{12}+z_4^3+z_5^3=0\bigr\}.
\end{equation*}
Calculation shows that $Y$ has trivial canonical bundle and three 
singular points $p_1=[0,0,0,0,1,-1]$, $p_2=[0,0,0,0,1,e^{\pi i/3}]$ 
and $p_3=[0,0,0,0,1,e^{-\pi i/3}]$, all modelled on~$\C^4/\an{i}$.

Now define a map $\si:Y\ra Y$ by
\begin{equation*}
\si:[z_0,\ldots,z_5]\longmapsto[\bar z_1,-\bar z_0,\bar z_3,
-\bar z_2,\bar z_5,\bar z_4].
\end{equation*}
Note that $\si^2=1$, though this is not immediately obvious,
because of the geometry of $\CP^5_{1,1,1,1,4,4}$. It can be shown 
that conditions (a)--(d) of Step 1 above hold for $Y$ and~$\si$.
\label{s3ex1}
\end{ex}

More suitable 4-folds $Y$ may be found by taking hypersurfaces
or complete intersections in other weighted projective spaces,
possibly also dividing by a finite group, and then doing a
crepant resolution to get rid of any singularities that we don't
want. Examples are given in \cite{Joyc4}, \cite[\S 15]{Joyc5}.

\subsection{Step 3: Resolving $\R^8/G$}
\label{s32}

Define $\al,\be:\R^8\ra\R^8$ by
\begin{equation*}
\begin{split}
\al:(x_1,\ldots,x_8)&\mapsto(-x_2,x_1,-x_4,x_3,-x_6,x_5,-x_8,x_7),\\
\be:(x_1,\ldots,x_8)&\mapsto(x_3,-x_4,-x_1,x_2,x_7,-x_8,-x_5,x_6).
\end{split}
\end{equation*}
Then $\al,\be$ preserve $\Om_0$ given in \eq{s1eq2}, so they lie in
$\Spin(7)$. Also $\al^4=\be^4=1$, $\al^2=\be^2$ and $\al\be=\be\al^3$.
Let $G=\an{\al,\be}$. Then $G$ is a finite nonabelian subgroup of
$\Spin(7)$ of order 8, which acts freely on $\R^8\sm\{0\}$. One
can show that if $Z$ is the compact $\Spin(7)$-orbifold constructed
in Step 2 above, then $T_{p_j}Z$ is isomorphic to $\R^8/G$ for
$j=1,\ldots,k$, with an isomorphism identifying the
$\Spin(7)$-structures $(\Om,g)$ on $Z$ and $(\Om_0,g_0)$ on $\R^8/G$,
such that $\be$ corresponds to the $\si$-action on~$Y$.

In the next two examples we shall construct two different ALE
$\Spin(7)$-manifolds $(X_1,\Om_1,g_1)$ and $(X_2,\Om_2,g_2)$
asymptotic to~$\R^8/G$. 

\begin{ex} Define complex coordinates $(z_1,\ldots,z_4)$ on $\R^8$ by
\begin{equation*}
(z_1,z_2,z_3,z_4)=(x_1+ix_2,x_3+ix_4,x_5+ix_6,x_7+ix_8),
\end{equation*}
Then $g_0=\ms{\d z_1}+\cdots+\ms{\d z_4}$, and $\Om_0=\ha\om_0
\w\om_0+\Re(\th_0)$, where $\om_0$ and $\th_0$ are the 
usual K\"ahler form and complex volume form on $\C^4$. In these 
coordinates, $\al$ and $\be$ are given by
\begin{equation}
\begin{split}
\al:(z_1,\ldots,z_4)&\mapsto(iz_1,iz_2,iz_3,iz_4),\\
\be:(z_1,\ldots,z_4)&\mapsto(\bar z_2,-\bar z_1,\bar z_4,-\bar z_3).
\end{split}
\label{s3eq1}
\end{equation}

Now $\C^4/\an{\al}$ is a complex singularity, as $\al\in\SU(4)$.
Let $(Y_1,\pi_1)$ be the blow-up of $\C^4/\an{\al}$ at 0. Then $Y_1$
is the unique crepant resolution of $\C^4/\an{\al}$. The action of
$\be$ on $\C^4/\an{\al}$ lifts to a {\it free} antiholomorphic
map $\be:Y_1\ra Y_1$ with $\be^2=1$. Define $X_1=Y_1/\an{\be}$. 
Then $X_1$ is a nonsingular 8-manifold, and the projection 
$\pi_1:Y_1\ra\C^4/\an{\al}$ pushes down to~$\pi_1:X_1\ra\R^8/G$. 

There exist ALE Calabi--Yau metrics $g_1$ on $Y_1$, which were
written down explicitly by Calabi \cite[p.~285]{Cal}, and are
invariant under the action of $\be$ on $Y_1$. Let $\om_1$ be
the K\"ahler form of $g_1$, and $\th_1=\pi_1^*(\th_0)$ the
holomorphic volume form on $Y_1$. Define $\Om_1=\ha\om_1\w
\om_1+\Re(\th_1)$. Then $(\Om_1,g_1)$ is a torsion-free 
$\Spin(7)$-structure on $Y_1$, as in Proposition~\ref{s1prop5}.

As $\be^*(\om_1)=-\om_1$ and $\be^*(\th_1)=\bar\th_1$, we 
see that $\be$ preserves $(\Om_1,g_1)$. Thus $(\Om_1,g_1)$ pushes 
down to a torsion-free $\Spin(7)$-structure $(\Om_1,g_1)$ on 
$X_1$. Then $(X_1,\Om_1,g_1)$ is an {\it ALE\/ $\Spin(7)$-manifold} 
asymptotic to~$\R^8/G$.
\label{s3ex2}
\end{ex}

\begin{ex} Define new complex coordinates $(w_1,\ldots,w_4)$ on $\R^8$ by
\begin{equation*}
(w_1,w_2,w_3,w_4)=(-x_1+ix_3,x_2+ix_4,-x_5+ix_7,x_6+ix_8).
\end{equation*}
Again we find that $g_0=\ms{\d w_1}+\cdots+\ms{\d w_4}$ and 
$\Om_0=\ha\om_0\w\om_0+\Re(\th_0)$. In these 
coordinates, $\al$ and $\be$ are given by
\begin{equation}
\begin{split}
\al:(w_1,\ldots,w_4)&\mapsto(\bar w_2,-\bar w_1,\bar w_4,-\bar w_3),\\
\be:(w_1,\ldots,w_4)&\mapsto(iw_1,iw_2,iw_3,iw_4).
\end{split}
\label{s3eq2}
\end{equation}
Observe that \eq{s3eq1} and \eq{s3eq2} are the same, except that 
the r\^oles of $\al,\be$ are reversed. Therefore we can use 
the ideas of Example \ref{s3ex2} again.

Let $Y_2$ be the crepant resolution of $\C^4/\an{\be}$. The action 
of $\al$ on $\C^4/\an{\be}$ lifts to a free antiholomorphic
involution of $Y_2$. Let $X_2=Y_2/\an{\al}$. Then $X_2$ is nonsingular, 
and carries a torsion-free $\Spin(7)$-structure $(\Om_2,g_2)$, making 
$(X_2,\Om_2,g_2)$ into an ALE $\Spin(7)$-manifold asymptotic to~$\R^8/G$.
\label{s3ex3}
\end{ex}

We can now explain the remarks on holonomy groups at the end of
Step 5. The holonomy groups $\Hol(g_i)$ of the metrics $g_1,g_2$
in Examples \ref{s3ex2} and \ref{s3ex3} are both isomorphic to
$\Z_2\lt\SU(4)$, a subgroup of $\Spin(7)$. However, they are two
{\it different} inclusions of $\Z_2\lt\SU(4)$ in $\Spin(7)$, as in
the first case the complex structure is $\al$ and in the second~$\be$.

The $\Spin(7)$-structure $(\Om,g)$ on $Z$ also has holonomy
$\Hol(g)=\Z_2\lt\SU(4)$. Under the natural identifications we have
$\Hol(g_1)=\Hol(g)$ but $\Hol(g_2)\ne\Hol(g)$ as subgroups of
$\Spin(7)$. Therefore, if we choose $i_j=1$ for all $j=1,\ldots,k$,
then $Z$ and $X_{i_j}$ all have the same holonomy group
$\Z_2\lt\SU(4)$, so they combine to give metrics $\ti g_t$ on
$M$ with~$\Hol(\ti g_t)=\Z_2\lt\SU(4)$.

However, if $i_j=2$ for some $j$ then the holonomy of $g$ on $Z$
and $g_{i_j}$ on $X_{i_j}$ are {\it different} $\Z_2\lt\SU(4)$
subgroups of $\Spin(7)$, which together generate the whole group
$\Spin(7)$. Thus they combine to give metrics $\ti g_t$ on $M$
with~$\Hol(\ti g_t)=\Spin(7)$.

\subsection{Conclusions}
\label{s33}

The author was able in \cite{Joyc4} and \cite[Ch.~15]{Joyc5} to
construct compact 8-manifolds with holonomy $\Spin(7)$ realizing
14 distinct sets of Betti numbers, which are given in Table
\ref{s7betti}. Probably there are many other examples which can
be produced by similar methods.

\begin{table}[htb]
\centering
{\caption{Betti numbers $(b^2,b^3,b^4)$ of compact
$\Spin(7)$-manifolds}\label{s7betti}}
{\begin{tabular}{ccccc}
\hline
\vphantom{$\bigr)^{k^k}$}
(4,\,33,\,200) & (3,\,33,\,202) & (2,\,33,\,204) & 
(1,\,33,\,206) & (0,\,33,\,208) \\ 
(1,\,0,\,908)  & (0,\,0,\,910)  & (1,\,0,\,1292) & 
(0,\,0,\,1294) & (1,\,0,\,2444) \\ 
(0,\,0,\,2446) & (0,\,6,\,3730) & (0,\,0,\,4750) & 
(0,\,0,\,11\,662)\vphantom{$\bigr)_{p_p}$} \\
\hline
\end{tabular}}
\end{table}

Comparing these Betti numbers with those of the compact 8-manifolds
constructed in \cite[Ch.~14]{Joyc5} by resolving torus orbifolds
$T^8/\Gamma$, we see that these examples the middle Betti number
$b^4$ is much bigger, as much as $11\,662$ in one case.

Given that the two constructions of compact 8-manifolds with holonomy 
$\Spin(7)$ that we know appear to produce sets of 8-manifolds with 
rather different `geography', it is tempting to speculate that the 
set of all compact 8-manifolds with holonomy $\Spin(7)$ may be rather
large, and that those constructed so far are a small sample with 
atypical behaviour.

\end{document}